\documentclass{amsart}
\usepackage[utf8]{inputenc}
\usepackage[margin=1.5in]{geometry}
\usepackage{hyperref}
\usepackage{url}
\usepackage{amsmath}
\usepackage{amsfonts}
\usepackage{amsthm}
\usepackage{amssymb}
\usepackage{xcolor}
\usepackage{graphicx}
\usepackage{mathtools}
\usepackage{tikz-cd}
\usepackage{ytableau}

\theoremstyle{plain}
\newtheorem*{ThA*}{Theorem A}
\newtheorem*{ThB*}{Theorem B}
\newtheorem*{ThC*}{Theorem C}
\newtheorem*{ThD*}{Theorem D}
\newtheorem*{ThE*}{Theorem E}
\newtheorem*{Con*}{Conjecture}
\newtheorem{Th}{Theorem}[section]

\newtheorem*{Cor*}{Corollary}
\newtheorem{Prop}[Th]{Proposition}
\newtheorem{Lemma}[Th]{Lemma}

\newtheorem{Ex}[Th]{Example}

\theoremstyle{definition}
\newtheorem{Def}[Th]{Definition}
\newtheorem{Asm}[Th]{Assumption}

\newtheorem{Rmk}[Th]{Remark}
\newtheorem{Set}[Th]{Setting}
\newtheorem{App}[Th]{Application}

\newcommand{\Hom}{\operatorname{Hom}}

\newcommand{\Spec}{\operatorname{Spec}}
\newcommand{\Proj}{\operatorname{Proj}}

\newcommand{\grade}{\operatorname{grade}}

\newcommand{\Tr}{\operatorname{Tr}}
\def\iff{if and only if }
\newcommand{\w}{\omega}

\def\ZZ{\mathbb Z}

\def\QQ{\mathbb Q}
\def\NN{\mathbb N}

\def\OO{\mathcal{O}}
\def\maxm{\mathfrak{m}}

\newcommand{\cO}{\mathcal{O}}

\newcommand{\fa}{\mathfrak{a}}

\newcommand{\maxM}{\mathfrak{M}}

\newcommand{\bbN}{\mathbb{N}}

\newcommand{\rS}{\mathrm{S}}

\newcommand{\bbZ}{\mathbb{Z}}

\newcommand{\bbQ}{\mathbb{Q}}

\newcommand{\bbR}{\mathbb{R}}

\newcommand{\la}{\longrightarrow}

\title{Test modules of extended Rees Algebras}
\author{Rahul Ajit, Hunter Simper}
\address{Department of Mathematics\\
	University of Utah\\
	Salt Lake City, UT 84112, USA.}
\email{rahulajit@math.utah.edu}

\begin{document}
\maketitle
\begin{abstract}
    Given a reduced, local ring $R$ and an ideal $\fa$ of positive height, we give a decomposition of the test module, $\tau(\w_T, t^{-\lambda})$, of the extended Rees algebra, $T =R[\fa t, t^{-1}]$. In particular, the degree zero component of this test module is $\tau(R, \fa^\lambda)$, thereby reducing the computation of test modules for non-principal ideals to the much easier case of principal ideals. Additionally, we apply our decomposition to generalize results on the $F$-rationality of Rees and extended Rees algebras \cite[Conjecture 4.1]{HWY-02F-Rational},\cite{Koley-Kummini21} as well as give simplified proofs of discreteness and rationality of F-jumping numbers, among other applications.

\end{abstract}

\section{Introduction}

Let $(R,\maxm)$ be a local ring of dimension $d\geq 2$ and $\mathfrak{a}\subseteq R$ be an ideal with positive height. Let $S=R[\mathfrak{a}t]$, $T=R[\mathfrak{a}t,t^{-1}]$, and $G = T/(t^{-1}) = \bigoplus_{n \geq 0}\fa^n/\fa^{n+1}$ be the Rees algebra, extended Rees algebra, and associated graded ring, respectively, of $R$ with respect to $\mathfrak{a}$. These rings carry deep geometric meaning, and relationships between various singularities, e.g., Cohen-Macaulayness and Gorensteinness, have been extensively studied see, \cite{GotoShimodaCM}, \cite{iai2024characterizationsgorensteinreesalgebras}, \cite{Goto-Nishida-Book}, \cite{Ikeda-Gor}, \cite{VietCM}, \cite{WhenCM}, \cite{WhenGor} \cite{Huneke-82} and, \cite{Lipman1994CM} among others. Of particular relevance to this paper, is $F$-rationality, which is measured by the test module, our main result being the following decompositions of the test modules of the pairs $(S, (\fa S)^\lambda)$ and $(T, (t^{-1})^\lambda)$:

\begin{ThA*}\label{Th-A}(see Theorem \ref{Th-test module decompositions} and, Theorem \ref{Th-testmoduleExtReesNeg})
Let $R$ be an F-finite reduced ring and $\fa \subseteq R$ be a ideal with height($\fa$) $> 0$. Let $S:=R[\mathfrak{a}t] =\oplus_{n \in \NN}{\mathfrak{a}}^nt^n$ denote the Rees algebra and $T := R[\mathfrak{a}t, t^{-1}] = \oplus_{n \in \ZZ} \mathfrak{a}^nt^n$ denote the extended Rees algebra where we define $\mathfrak{a}^n := R$, for $n \leq 0$ and $t$ be a dummy variable of degree $1$, and $\lambda \geq 0$.
     Then, the test submodule of $(S, (\fa S)^\lambda)$ is given by
\[\tau(\omega_{S}, (\mathfrak{a}.S)^{\lambda}) = \bigoplus_{n\geq 1} \tau (\omega_R, \mathfrak{a}^{n+\lambda})\] and the test submodule of $(T, (t^{-1})^\lambda)$  is given by
\[ [\tau(\omega_{T}, (t^{-1})^\lambda)]_{\geq 0} = \bigoplus_{n \geq 0} \tau (\omega_R, \mathfrak{a}^{n+\lambda}).\]
For $n \leq -1$,
    \[[\tau(\omega_T, (t^{-1})^\lambda)]_n = 
    \begin{cases} 
\tau(\omega_R) , & \forall n \leq -\lambda \\
\tau(\omega_X, \mathfrak{a}^{\lambda+n}), & -\lambda < n \leq -1
\end{cases}.
    \]
 In particular, $ [\tau(\omega_{T}, (t^{-1})^\lambda)]_{\geq 1} = \tau(\omega_{S}, (\mathfrak{a}.S)^{\lambda}) $ and $[\tau(\omega_{T}, (t^{-1})^\lambda)]_0 = \tau(\omega_R , \mathfrak{a}^\lambda)$.
\end{ThA*}

Similar decompositions of test modules for the Rees algebra $S$ have been computed in the case of strong assumptions on $R$, $S$, and $\mathfrak{a}$, see \cite[Theorem 5.1]{HaraYoshida}, \cite[Theorem 1.1]{KotalKummini}. To the authors' knowledge, no decomposition of test modules $T$ has been previously described. Analogous descriptions of multiplier modules have been computed (also in restrictive settings) in \cite[Proposition 3.1]{HyryBlowUp}, and in forthcoming work, the first author provides a characteristic 0 and mixed characteristic analogue to Theorem A, see \cite{AjitChar0Mixed}, (also see \cite[Theorem 1 and (1.3.1)]{BMS-06-BS}, for a similar decomposition in characteristic 0).


The decomposition of Theorem A has both computational and theoretical applications. On the computational front, while there exists algorithms for computing test modules for principal ideals prior to this decomposition, there was not an method to verifiably compute a test module for a non-principal ideal. This will be implemented in a forthcoming package for Macaulay 2 \cite{M2}.

On the theory side, when $\lambda=0$, this decomposition allows the direct comparison on the test modules and hence the $F$-rationality of $S$ and $T$. While similar questions about the rationality of $S$ and $T$ have been studied by many authors, c.f., \cite{HWY-F-Regular}, \cite{HyryCoef} and \cite{LipmanAdjoint}, the work on $F$-rationality has been restricted to a $\maxm$ primary ideal $\mathfrak{a}$. In \cite{HWY-02F-Rational} Hara, Watanabe, and Yoshida proved that for $\maxm$-primary ideal $\fa$ and $F$-rational ring R, $S$ is $F$-rational implies $T$ is $F$-rational and posed the converse as a conjecture, stated fully, their conjecture is as follows.

\begin{Con*} (\cite[Conjecture 4.1]{HWY-02F-Rational}, now a theorem of \cite{Koley-Kummini21})
    Let $R$ be an excellent $F$-rational local ring and $\fa \subseteq R$ be an $\maxm$-primary ideal of $R$. Then $T=R[\mathfrak{a}t, t^{-1}]$ is F-rational implies $S=R[\mathfrak{a}t]$ is F-rational.
\end{Con*}

A few words about this conjecture: If Boutot-type theorem (which is true for rational singularities in characteristic 0, see \cite{Boutot}) would hold for $F$-rational singularities, then this conjecture would follow trivially, see \cite[Page 181]{HWY-02F-Rational}. However, Boutot-type theorem fails in positive characteristics, see \cite{Watanabe-Boutot}. Geometrically, Spec $T$ corresponds to the deformation to the normal cone Spec $G$, see \cite[Chapter 5.]{Fulton}. A positive answer to this conjecture implies the following: if $G$ is $F$-rational then so is $S$, see \cite[see page 182]{HWY-02F-Rational} and Remark \ref{deformation}.

With the stated assumptions, the conjecture was first proved by \cite{Koley-Kummini21}, using different techniques. We will give a positive answer to (a generalization of) this conjecture when $\fa$ is not necessarily $\maxm$-primary, see Theorem \ref{Th-equivofFrationality}, as well as a new proof of \cite[Theorem 4.2]{HWY-02F-Rational} with this relaxed assumption. 

\begin{ThB*}\label{Th-B}(see Theorem \ref{Th-equivofFrationality})
 With the notations introduced above, we have the following:
    \begin{enumerate}
        \item $R$ and $S$ are both $F$-rational implies $T$ is $F$-rational.
        \item $T$ is $F$-rational implies $R$ and $S$ are $F$-rational.
    \end{enumerate} 
\end{ThB*}

In order to prove both Theorems A and B, it is necessary to compute and compare canonical modules $\w_\rS$ and $\w_T$. Such comparisons were previously known when $\fa$ is $\maxm$-primary \cite[Proposition B.6]{MSTWW}, or when height$(\fa) \geq 2$ \cite[Corollary 3.3]{TomariWatanabe} and we extend this comparison to its fullest generality.


\begin{ThC*}\label{Th-C} (see Theorem \ref{Th-compareomegaSandT}, Proposition \ref{Prop-negativeeonghdegrees} and Proposition \ref{Prop-canonicalnegativedegrees})
   With the above notation, assume $\omega_S$ and $\omega_T$ exist, e.g., when $R$ is $F$-finite \cite{Gabber}, then 
   $[\omega_S]_n=[\omega_T]_n$ for all $n>0$. For all $n \ll 0$ we have that $[\omega_T]_n \simeq \omega_R$. When $S$ is Cohen-Macaulay, $[\omega_T]_i \simeq \omega_R$ for all $i\leq 0$ and $\omega_T$ is generated in non-negative degrees.
\end{ThC*}


As additional applications, in Section \ref{sec-other applications}, we highlight how Rees' ``principalization trick" in the form of Theorem A allows for reduction of problems to the principal ideal case and provide an alternate proof of the discreteness and rationality of $F$-jumping numbers, \cite[Theorem B]{Schwede-Tucker2014Testidealsofnon-principalideals:computationsjumpingnumbersalterationsanddivisiontheorems}.


\section{Acknowledgments}

First and foremost, we would like to thank Karl Schwede who organized the RTG Macaulay2 semester and suggested the project of computing test ideal for non-principal ideals in Macaulay2, and provided many helpful clarifications. The first-named author would like to thank his advisors, Christopher Hacon and Karl Schwede, for their constant encouragement, unwavering support, inspiring teachings, and infinite patience. He also thanks Daniel Apsley, Harold Blum, Brad Dirks, David Eisenbud, Linquan Ma, Mircea Mustaţă, and Kevin Tucker for encouraging discussions while presenting this paper. Finally, he would like to thank 
Joaqu\'in Moraga for carefully reading this manuscript and providing many valuable feedback and helpful comments, which improved the exposition substantially, and Bernd Ulrich for kindly sharing his beautiful notes on Local Cohomology.

\section{Canonical Modules}

\begin{Asm}\label{Assumption}
Throughout this section let $(R,\maxm)$ be a local ring of dimension $d\geq 2$ and $\mathfrak{a}\subseteq R$ an ideal with positive height. We set $\rS=R[\mathfrak{a}t]$ to be the Rees algebra of $\mathfrak{a}$ with homogeneous maximal ideal $\maxm_S =  \maxm \oplus \mathfrak{a} \oplus \mathfrak{a}^2 \oplus \dots$ and $T=R[\mathfrak{a}t,t^{-1}]$ the extended Rees algebra of $\mathfrak{a}$ with homogeneous maximal ideal $\maxm_T$.
\end{Asm}
\begin{Def}
Given a ring $R$ and an ideal $J \subseteq R$, the local cohomology with respect to $J$ is the set of derived functors \(H^{i}_{J}(-)\) on the category of \(R\)-modules for the functor \(H^{0}_{J}(-)\) defined by

\[H^{0}_{J}(M)=\{m\in M \mid m|_{{\mathcal{U}}}=0\ \mbox{ as a section of }\widetilde{M},\mbox{ where }\ {\mathcal{U}}=\operatorname{Spec}R\setminus\mathbb{V}(J)\},\]

\(\widetilde{M}\) denotes the quasicoherent sheaf on \(\operatorname{Spec}R\) determined by \(M\). If \(J\) is a finitely generated ideal of \(R\), we can also write

\[H^{0}_{J}(M)=\bigcup_{t\in\mathbb{N}}\{m\in M \mid J^{t}m=0\}.\]
\end{Def}

\begin{Rmk}\label{RemarkSdS}

The Sancho de Salas sequence for the Rees Algebra $S$ \cite[Section 2.5.2]{HyrySmith} is 
$$
\dots \la H^i_{\maxm_\rS}(\rS) \la  \bigoplus_{n\in \bbN} 
H^i_{\maxm}(\fa^n) \la 
\bigoplus_{n \in \bbZ} H^i_{Z}(Y, \cO_Z(n)) \la
H^{i+1}_{\maxm_\rS}(\rS) \la \dots,
$$ where $\pi: Y = \Proj \rS \to X$ the blowup along $\mathfrak{a}$ so that $\mathfrak{a} \OO_Y= \OO_Y(-G)$ and $Z = Y \times_{\Spec R} \Spec (R/\maxm)$ is the scheme-theoretic fiber over the closed point $\maxm$ of 
$\Spec A$.
We get, as graded $A$-modules, 
$H^{d+1}_{\maxm_\rS}(\rS) \cong \oplus_{n < 0}H^d_Z(Y, \cO_Y(-nG))$ \cite[2.5.2 (1)]{HyrySmith} . This is because 
the maps $H^d_m(\fa^n) \la H_Z^d(Y, \cO_Y(-nG))$ are surjective for all
 $n \geq 0 $
(see, \cite[page 103]{LipmanTeissier}). Hence, Matlis dually (see \ref{defGradedCanonical}) gives, \[\omega_\rS \simeq\bigoplus_{n>0}H^0(Y,\omega_Y(-nG)).\]

\end{Rmk}

\begin{Def}\cite[Definition 3.1]{Kim}, \label{defGradedCanonical}
    Let $A$ be a Noetherian, $\ZZ$-graded ring of dimension d+1 and $\maxM$ be the unique maximal homogeneous ideal. Let $E =E(A_0/\maxm)$ be the injective hull of the residue field of the local subring $(A_0, \maxm)$. For $R$-modules $M$ and $N$, define $^*\mathrm{Hom}_A(M , N)$ to be the R-submodule of $\mathrm{Hom}_R(M, N)$ generated by the homogeneous $R$-linear maps of arbitrary degree and $\widehat{M} := M \otimes_{A_0} \widehat{A}^\maxm$. Note, for finitely generated modules $\widehat{M} \cong \widehat{N}$ implies $M \cong N.$ Finally, we say $\w_A$ is the graded canonical module of A if $\widehat{\w_A}\cong$  $^*\mathrm{Hom}_{A_0}(H^{d+1}_\mathfrak{M}(A), E)$
\end{Def}

\subsection{Comparing $\omega_S$ and $\omega_T$}
The objective of this section is to compare the graded canonical modules of $S$ and $T$. To this end our main result is that $\omega_S$ and $\omega_T$ agree in all positive degrees and we will spend the remainder of the section proving our main theorem.

\begin{Th}\label{Th-compareomegaSandT}
Under Assumption \ref{Assumption}, $[\omega_S]_n=[\omega_T]_n$ for all $n>0$. 
\end{Th}

To compare the graded structure of $\omega_S$ and $\omega_T$ we will analyze the graded structure of $H^{d+1}_{\maxm_\rS}(\rS)$ and $H^{d+1}_{\maxm_T}(T)$. To do this we will make use of the following long exact sequences:

\begin{equation}\label{LES-quot}
    \begin{tikzcd}
        \cdots \ar{r} &H^d_{\maxm_\rS}(\rS) \ar{r} & H^d_{\maxm_\rS}(T) \ar{r} & H^d_{\maxm_\rS}(\rS/T) \ar{r} 
        & H^{d+1}_{\maxm_\rS}(\rS) \ar{r} & H^{d+1}_{\maxm_\rS}(T) \ar{r} & \cdots
    \end{tikzcd}
\end{equation}
and 
\begin{equation}\label{LES-add}
    \begin{tikzcd}
        \cdots \ar{r}& H^d_{\maxm_T}(T) \ar{r}  &H^d_{\maxm_\rS}(T) \ar{r} & H^d_{{\maxm_\rS}_{t^{-1}}}(T_{t^{-1}}) \ar{r} & H^{d+1}_{\maxm_T}(T) \ar{r} & H^{d+1}_{\maxm_\rS}(T) \ar{r} & \cdots
    \end{tikzcd}
\end{equation}
where (\ref{LES-quot}) is induced by the short exact sequence of $\rS$-modules \[\begin{tikzcd}
    0 \ar{r} & \rS \ar{r} & T \ar{r} & \rS/T \ar{r} & 0
\end{tikzcd}\] 
and (\ref{LES-add}) is induced from the short exact sequence of complexes
\[\begin{tikzcd}
        0 \ar{r} & \check{C}^\bullet([\underline{g}];T)\otimes T_{t^{-1}} [1] \ar{r} & \check{C}^\bullet([\underline{g},t^{-1}];T) \ar{r}& \check{C}^\bullet([\underline{g}];T) \ar{r} & 0
    \end{tikzcd} \]
    where we take $\underline{g}$ to be a generating set for $\maxm_\rS$ and note that $\maxm_T=\maxm_\rS+(t^{-1})$. 


    We will need to following lemmas describing the vanishing and $R$-module structure of the terms in these sequences.

\begin{Lemma}\label{Lemma-localcohomologiesISO1}
    \begin{enumerate}
        \item  For all $k$ we have an isomorphism:
        \[H^k_{\maxm_\rS}(T_{t^{-1}}) \simeq \bigoplus_{i\in \ZZ} H^k_\maxm(R) t^i.\]
        In particular $H^{d+1}_{\maxm_\rS}(T_{t^{-1}})=0$.

        \item For $k\geq 2$ we have an isomorphism:
    \[H^k_{\maxm_\rS}(T/\rS) \cong \bigoplus_{i<0} H^k_\maxm(R) t^i.\]
    \end{enumerate}
\end{Lemma}

\begin{proof}
    \begin{enumerate}
        \item  The claim follows from the fact that $T_{t^{-1}}= R[t,t^{-1}]$ and ${\maxm_S} T_{t^{-1}}=\maxm_S R[t,t^{-1}]=\maxm R[t,t^{-1}]$.

        \item Let $S_+=(\fa t) \subseteq S$ be the ideal of elements of $S$ of positive degree. Then $\maxm_S=\maxm+S_+$ so we have the Mayer-Vietoris sequence:
    \[\begin{tikzcd}
        H^{k-1}_{\maxm S\cap S_+}(T/S) \ar{r} & H^{k}_{\maxm_S}(T/S) \ar{r} & H^{k}_{\maxm}(T/S) \oplus H^{k}_{S_+}(T/S) \ar{r} & H^{k}_{\maxm S\cap S_+}(T/S) 
    \end{tikzcd}\]
    Note that $T/S$ is a $S_+$-torsion module and a $\maxm\cap S_+$-torsion module, hence since $k$ is at least 2 so $H^{k-1}_{\maxm S\cap S_+}(T/S)$, $H^{k}_{S_+}(T/S)$ and $H^{k}_{\maxm S\cap S_+}(T/S)$ are all zero. Thus we have that
    \[H^{d}_{\maxm_S}(T/S) \cong H^{d}_{\maxm}(T/S) \cong  H^d_\maxm(\bigoplus_{i<0} R t^i ) \cong \bigoplus_{i<0} H^d_\maxm(R) t^i. \]
    \end{enumerate}
\end{proof}

\begin{Lemma}\label{Lemma-localcohomologiesISO2}
   With notation as above,
   \[H_{\maxm_S}^{d+1}(T)=0.\]
\end{Lemma}

\begin{proof}
We prove that the ideal $\maxm_S T$ is not $\maxm_T$ primary and so the claim follows from Hartshorne-Lichtenbaum Vanishing, \cite[Tag 0EB0]{stacks-project}. Observe that for $i\leq 0$ the degree $i$ component of $\maxm_S T$ is isomorphic to $\maxm$, hence in each non-positive degree $T/ \maxm_S T$ is isomorphic to $R/\maxm$. Since this is true for all non-positive degrees we have that $T/ \maxm_S T$ is not a finite dimensional vector space.
%
%
%
\end{proof}

We are now ready for the main technical result of this section, which will allow us to prove Theorem \ref{Th-compareomegaSandT}.

\begin{Prop}\label{Prop-5-lemma-setup}
    There is a commutative diagram of $R$-modules:
\begin{equation}   
\label{localcohomologydiagram}
\begin{tikzcd} 
H_{\maxm_S}^{d}(T) \ar["\alpha"]{r} \ar["\phi_1"]{d}& \oplus_{i<0} H_{\maxm}^{d}(R) t^i \ar["\beta"]{r} \ar["\phi_2"]{d} & H_{\maxm_S}^{d+1}(S) \ar{r} \ar["\phi_3"]{d} & 0\\
H_{\maxm_S}^{d}(T) \ar["\alpha' "]{r} & \oplus_{i\in \ZZ} H_{\maxm}^{d}(R) t^i \ar["\beta' "]{r} & H_{\maxm_T}^{d+1}(T) \ar{r} & 0,
\end{tikzcd}
\end{equation} 
 where the top row comes from \eqref{LES-quot}, the bottom row comes from \eqref{LES-add} and the vertical maps are as follows:
 \begin{itemize}
     \item $\phi_1$ is equality in negative degrees and zero otherwise
     \item $\phi_2$ is the natural inclusion
     \item $\phi_3$ is induced by the commutativity of the left square.
 \end{itemize}
 \end{Prop}

\begin{proof}
The exactness of rows follows from the exactness of \eqref{LES-quot} and  \eqref{LES-add} along with Lemma \ref{Lemma-localcohomologiesISO1} and Lemma \ref{Lemma-localcohomologiesISO2}. By construction, to show the diagram commutes it is sufficient for show that left square commutes.

Fix a generating set of $ \mathfrak{a}=(\underline{f})$ and generating set $\underline{x}$ of $\maxm$. Then $\maxm_S=(\underline{ft})+(\underline{x})$ so an element of $H^d_{\maxm_S}T$ of degree $n$ is of the form:
\[\left[\sum_\sigma \frac{f_\sigma t^{i_\sigma}}{x^{a_\sigma} f^{b_\sigma}t^{|b_\sigma|}}\right]\]
where $i_\sigma+|b_\sigma|=n$, and $\sigma$ runs over the number of generator of $\fa$.

Then since $\alpha$ is the natural map $H^d_{\maxm_S}(T) \to H^d_{\maxm_S}(S/T)$ we see that
\[ \alpha\left(\left[\sum_\sigma \frac{f_\sigma t^{i_\sigma}}{x^{a_\sigma} f^{b_\sigma}t^{|b_\sigma|}}\right] \right)=\begin{cases}
    0 & \text{ if } i_\sigma+|b_\sigma|=n\geq 0\\
    \left[\sum_\sigma \frac{f_\sigma }{x^{a_\sigma} f^{b_\sigma}}\right] t^{n} & \text{ if } i_\sigma+|b_\sigma|=n< 0
\end{cases}.
\]
On the other hand, $\alpha'$ is just the localization at $t^{-1}$ map so
\[\phi_1\left( \left[\sum_\sigma \frac{f_\sigma t^{i_\sigma}}{x^{a_\sigma} f^{b_\sigma}t^{|b_\sigma|}}\right]\right)=\begin{cases}
    0 & \text{ if } i_\sigma+|b_\sigma|=n\geq 0\\
    \left[\sum_\sigma \frac{f_\sigma t^{i_\sigma}}{x^{a_\sigma} f^{b_\sigma}t^{|b_\sigma|}}\right] & \text{ if } i_\sigma+|b_\sigma|=n< 0
\end{cases} 
\]
So 
\[\alpha'\left(\phi_1\left(\left[\sum_\sigma \frac{f_\sigma t^{i_\sigma}}{x^{a_\sigma} f^{b_\sigma}t^{|b_\sigma|}}\right]\right)\right)=\begin{cases}
    0 & \text{ if } i_\sigma+|b_\sigma|=n\geq 0\\
    \left[\sum_\sigma \frac{f_\sigma }{x^{a_\sigma} f^{b_\sigma}}\right] t^{n} & \text{ if } i_\sigma+|b_\sigma|=n< 0
    \end{cases}.
\]
Hence, the left square commutes.
    
\end{proof}

Theorem \ref{Th-compareomegaSandT} now follows immediately from Proposition\ref{Prop-5-lemma-setup}.

\begin{proof}[Proof of Theorem \ref{Th-compareomegaSandT}]
    Restricting the diagram \eqref{localcohomologydiagram} of Proposition \ref{Prop-5-lemma-setup} to negative degrees we have that that $[H^{d+1}_{\maxm_S}(S)]_n\simeq [H^{d+1}_{\maxm_T}(T)]_n$ for all $n<0$. The claim then follows by dualizing.
\end{proof} 

\subsection{Negative degrees and generators of $\omega_T$}
In general it is hard to describe negative degrees of $\omega_T$, however asymptotically it is well understood.

\begin{Prop}\label{Prop-negativeeonghdegrees}
For all $n \ll 0$ we have that $[\omega_T]_n \simeq \omega_R$.
\end{Prop}
\begin{proof}
    Local cohomology and hence canonical modules localize, thus we have an inclusion $\omega_T \subseteq \omega_{T_{t^{-1}}} = \omega_{R[t,t^{-1}]} = \bigoplus_{n\in \ZZ} \omega_R t^n$. Fix a generating set $x_1,\ldots, x_l$ of $\omega_Rt^0 = \w_R$ viewed as elements of $[\omega_{T_{t^{-1}}}]_0$. By $\omega_T \subseteq \omega_{T_{t^{-1}}} = \omega_{R[t,t^{-1}]} = \bigoplus_{n\in \ZZ} \omega_R t^n$, we see that for each $i$ there exists some non-negative integer $k_i$ and $y_i \in [\w_T]_{-k_i}$ such that $x_i = y_it^{k_i}$, giving $x_i (t^{-1})^{k_i} \in [\omega_T]_{-k_i}$. Thus for all $k\geq \max(k_1,\ldots, k_l)$, $x_i t^{-k} \in [\omega_T]_{-k}$ and hence $\omega_R \subseteq [\omega_T]_{-k} \subseteq \omega_R$, giving $[\omega_T]_{-k} = \omega_R$.
\end{proof}

\begin{Rmk}\label{inclusionofgradedpieces}
    As $\omega_T \subseteq \omega_{T_{t^{-1}}} = \omega_{R[t,t^{-1}]} = \bigoplus_{n\in \ZZ} \omega_R t^n$, we can write $\omega_T = \bigoplus_{n\in \ZZ} \omega_nt^n$ where $\w_n \subseteq \w_R$ is an $R$-submodule. Note, $\forall n \in \ZZ$, multiplication by $t^{-1}$ is just the inclusion map $\w_n \subseteq \w_{n-1}\subseteq\w_R$. Hence, $[\w_T \cdot t^{-N}]_k = [\w_T]_{N+k}$. We will use this simple trick many times in this paper.
\end{Rmk}

When $S$ is Cohen-Macaulay and we can use diagram \eqref{localcohomologydiagram} in Proposition \ref{Prop-5-lemma-setup} (augmented on the left) to finely control where this stability of negative components of $\omega_T$ occurs as well as describe its generating degrees. 

\begin{Prop}\label{Prop-canonicalnegativedegrees}
    Suppose $S$ is Cohen-Macaulay. Then, $[\omega_T]_i \simeq \omega_R$ for all $i\leq 0$ and $\omega_T$ is generated in non-negative degrees.
\end{Prop}
\begin{proof}
    Consider the diagram \eqref{localcohomologydiagram} of Proposition \ref{Prop-5-lemma-setup}. Since $S$ is Cohen-Macaulay of dimension $d+1$ we have that the map $\alpha$ is injective, hence $H_{\maxm_S}^d(T)$ is concentrated in negative degrees. The map $\alpha'$ is degree $0$ hence the image is concentrated in negative degrees. Thus $\beta'$ is an isomorphism in non-negative degrees and the first part of the proposition follows after applying Graded Matlis duality. For the statement on generating degrees, we note that from the above computation the highest possible degree of a socle element of $H_{\maxm_S}^d(T)$ is degree 0.
\end{proof}

\begin{Rmk}
    Similar comparisons of canonical modules of $S$ and $T$ under restricted setups appeared in the past. For example, see \cite[Corollary 3.3]{TomariWatanabe} when $height(\fa) \geq 2$ and \cite[Appendix B]{MSTWW} when $\fa$ is $\maxm$-primary.
\end{Rmk}

\begin{Rmk}
    We would like to mention that all results proved in this section are characteristic-free.
\end{Rmk}

\section{Positive Characteristic preliminaries}

In this section, we collect some preliminaries, all of which are well-known to experts. For further details, we refer to \cite{FBook}. 

\begin{Set}

Let \(R\) be a reduced Noetherian ring of characteristic \(p > 0\) with perfect residue field. The \textit{absolute Frobenius} \(F: R \to R\) is the ring homomorphism defined by \(F(r) = r^p\), which we assume to be finite. For \(R\)-modules, we define the \textit{Frobenius pushforward} \(F_* M\) which is \(M\) as an abelian group with twisted \(R\)-action: \(r \cdot m = r^p m\). X  will denote Spec R.
\end{Set}

\begin{Rmk}
    A local ring $(R, \maxm, k)$ is $F$-finite \iff $R$ is excellent and $k$ is $F$-finite \cite{Kunz}.
\end{Rmk}

\begin{Rmk}
    When $R$ is $F$-finite, then by a result of Gabber \cite{Gabber}, $R$ is a
homomorphic image of an F-finite regular ring and hence admits a canonical module.
\end{Rmk}

\begin{Def}\label{Def-trace}
    The Grothendieck dual to the Frobenius map $\cO_X \to F^e_* \cO_X$ gives a map $F^e_* \omega_X^{\bullet} \to \omega_X^{\bullet}$ on the (normalized) dualizing complex.  By taking the $(-d)th$ cohomology, we get the trace map on canonical modules $Tr_F$:$F^e_* \omega_X \cong Hom_{\cO_X}(F^e_*\cO_X, \w_X) \to \omega_X$, which can be thought of as evaluation-at-1 map.
\end{Def}

\begin{Rmk}
    For an $F$-finite domain $R$, $Tr_F$ generates $\Hom_{R}(F^e_* \omega_R, \omega_R)$ as an $F^e_*R$-module  : $\Hom_{R}(F^e_* \omega_R, \omega_R) = \langle \Tr^e_F\rangle_{ F^e_* R} \cong F^e_* R$
\end{Rmk}

\begin{Ex}
    [Trace for Polynomial Rings]
Let \(R = k[x_1, \dots, x_n]\) with perfect \(k\), \(char(k) = p\). The trace map \(\Tr: R^{1/p} \to R\) is explicitly given by:
\[
\Tr\left( x_1^{a_1/p} \cdots x_n^{a_n/p} \right) = 
\begin{cases} 
x_1^{(a_1 - (p-1))/p} \cdots x_n^{(a_n - (p-1))/p} & \text{if } p \mid (a_i + 1) \ \forall i \\
0 & \text{otherwise}
\end{cases}
\]
for exponents \(a_i \in \mathbb{Z}\).

\end{Ex}

\begin{Prop}\label{Prop-traceonrees}

\begin{enumerate}
    \item The map trace on $\Tr_T: [\omega_T]_{n p} \to [\omega_T]_n$ is the restriction of the trace map $\Tr_R: \omega_R \to \omega_R$.
    \item If $\grade(\mathfrak{a})>0$ then the map trace on $\Tr_S: [\omega_S]_{n p} \to [\omega_S]_n$ is the restriction of the trace map $\Tr_R: \omega_R \to \omega_R$
\end{enumerate}
    .
\end{Prop}
\begin{proof}
    For the claim of $\Tr_T$: First consider the inclusion $\omega_T \subseteq \omega_{T_{t^{-1}}} = \omega_{R[t,t^{-1}]} = \bigoplus_{n\in \ZZ} \omega_R t^n$ gives us that  $[\omega_T]_n \subseteq \omega_R$ and note that $\Tr_{T_{t^{-1}}}$ clearly restricts to $\Tr_R$ on graded components. The claim then follows from the fact that $\Tr_T$ generates $\Hom_T(F_* \omega_T, \omega_T)$, and hence $(\Tr_T)_{t^{-1})}$ generates $(\Hom_T(F_* \omega_T, \omega_T))_{t_{-1}} \simeq \Hom_{T_{t^{-1}}}(F_* \omega_{T_{t^{-1}}}, \omega_{T_{t^{-1}}})$, i.e., $(\Tr_T)_{t^{-1})} = \Tr_{T_{t^{-1}}}$.
    
    For claim on $\Tr_S$: let $a\in \mathfrak{a}$ be a nonzero-divisor, then $S_{at} = R[t,t^{-1}]$ and the argument follows mutatis mutandis. 
\end{proof}

\begin{Rmk}
     If $R$ the Frobenius action on $H^d_\maxm(R)$ injective, e.g. when $R$ is $F$-injective, then $[\omega_T]_n = \omega_R$ for all $n < 0$.
\end{Rmk}

\begin{Def}
    A normal Cohen-Macaulay local ring \((R, \maxm)\) of dimension \(d\) has \textit{pseudo-rational singularities} if for every proper birational morphism \(\pi: X \to \Spec R\) with \(X\) normal, the natural map:
\[
H^d_{\maxm}(R) \to H^d_{\pi^{-1}(\maxm)}( \mathcal{O}_X)
\]
is injective.

\end{Def}

\begin{Def}
    An F-finite local ring \((R,\mathfrak{m})\) of prime characteristic is \textit{\(F\)-injective} if the natural Frobenius action on local cohomology \(H_{\mathfrak{m}}^{i}(R)\xrightarrow{F^{e}}F_{*}^{e}H_{\mathfrak{m}}^{i}(R)\) is injective for all \(i\in\mathbb{N}\) and all \(e>0\).

    A Noetherian local ring \((R,\mathfrak{m})\) of prime characteristic and dimension \(d\) is \textit{\(F\)-rational} if \(R\) is Cohen-Macaulay and for all non-zerodivisors \(c\), the map on local cohomology

\begin{equation}\label{F-loc}
H_{\mathfrak{m}}^{d}(R)\xrightarrow{F_{*}^{e}c\circ F^{e}}F_{*}^{e}H_{\mathfrak{m}}^{d}(R)
\end{equation}

is injective for some \(e>0\) (equivalently, for all \(e\gg 0\)).
    
\end{Def}

\begin{Def}
    A ring \(R\) of characteristic \(p > 0\) is \textit{strongly \(F\)-regular} if for every \(c \in R^\circ\), there exists \(e \geq 1\) such that the composition:
\[
R \xrightarrow{\cdot c} R \xrightarrow{F^e} F^e_* R
\]
splits as \(R\)-module homomorphisms.
\end{Def}
\begin{Rmk}
    Strongly \(F\)-regular rings are normal, Cohen-Macaulay and \(F\)-rational (hence pseudo-rational). Moreover, any Gorenstein $F$-rational ring is strongly F-regular.
\end{Rmk}

\begin{Def}
    A strong test element for $R$ is an element c with the property that for all non-zerodivisor d, there exists an $e_0 >0$  such that c is in the image of the map
    $$Hom_{R}(F^e_*R , R)\xrightarrow{\textit{evaluation at $F^{e}_*d$}}R$$ for all $e \geq e_0$. Note, $R$ is strongly F-regular \iff 1 is a strong test element. For a Noetherian F-finite reduced ring, it is an important result of Hochster and Huneke that a strong test element which is also a non-zero divisor, exists.
\end{Def}

\begin{Def}
    By a triple $(X, \Delta, \fa^t)$, $t \in \bbR_{\geq 0}$ we mean a normal integral scheme $X$ together with an effective $\QQ$-divisor $\Delta$ on $X$ such that $K_X + \Delta$ is $\bbQ$-Cartier, and an ideal $\fa \subseteq \cO_X$.
\end{Def}

\begin{Def}\label{Def-testideal}

    Let $(X, \Delta, \fa^\lambda)$ is a triple where $X = \Spec R$. The test ideal $\tau(X, \Delta, \fa^\lambda)$ is the unique smallest non-zero ideal $J \subseteq R$ such that for every $e > 0$ and every section $\phi \in \Hom_{R}(F^e_* R(\lceil (p^e - 1)\Delta \rceil), R) \supseteq \Hom_{R}(F^e_* R, R)$, we have $\phi(F^e_* (\fa^{\lceil t(p^e - 1) \rceil} J)) \subseteq J$.

\end{Def}

\begin{Def}\label{Def-Testmodule}
    Assume the setup and let $\Gamma \geq 0$ be a $\bbQ$-divisor. Then the test module, $\tau(\omega_R, \Gamma, \fa^{t})$, is the unique smallest nonzero submodule $J \subseteq \omega_X$ such that
for every $e \geq 0$ and every map $\phi \in \Hom_{R}\big(F^e_* (\omega_R( \lceil (p^e - 1) \Gamma \rceil)), \omega_R\big)$
 we have that $\phi\big(F^e_* ( \fa^{\lceil \lambda(p^e - 1) \rceil}\cdot J)\big) \subseteq J.$
When $\Gamma = 0$ (or $\fa = R$) we write, $\tau(\omega_R,\fa^{\lambda})$ (or $\tau(\omega_R, \Gamma)$, respectively).
\end{Def}

\begin{Rmk}
    It is an easy exercise \cite[Chapter-2, Section 5.2]{FBook} to see that for an F-finite reduced ring $R$ with canonical module $\w_R$ and a strong test element c which is also a non-zero divisor, $\tau(\w_R) = Tr^e_R(F^e_*(c.\w_R))$ for $e\gg0.$
\end{Rmk}

\begin{Rmk}
    We also note that a reduced F-finite ring $R$ is $F$-rational \iff $R$ is Cohen-Macaulay and $\tau(\w_R) = \w_R.$
\end{Rmk}

\begin{Lemma} \cite[Lemma 4.2]{Schwede-Tucker2014Testidealsofnon-principalideals:computationsjumpingnumbersalterationsanddivisiontheorems} \label{Relatingtestmoduleandideal} 
   Choose $K_R$ such that $- K_R$ is effective and fix $\omega_R = \cO_X(K_R) \subseteq K(X)$, then
\[
\tau(\omega_R, \Gamma, \fa^{\lambda}) = \tau(R, \Gamma-K_R, \fa^{\lambda})
\]
or, equivalently if $\Delta = \Gamma -K_R$, then
\[
\tau(\omega_R, \Delta + K_R, \fa^{\lambda}) = \tau(R, \Delta, \fa^{\lambda}).
\]
\end{Lemma}

Following main result from \cite{Schwede-Tucker2014Testidealsofnon-principalideals:computationsjumpingnumbersalterationsanddivisiontheorems} will be the most crucial tool for us.

\begin{Prop}\cite[Proposition 3.6]{Schwede-Tucker2014Testidealsofnon-principalideals:computationsjumpingnumbersalterationsanddivisiontheorems}\label{Prop-tau from blowup} 
Let $\mathfrak{a}$ be any ideal of R. Consider $\pi : Y \to X$, the normalized blowup of $X$ along $\mathfrak{a}$, so that $\mathfrak{a} \OO_{Y} = \OO_{Y}(-G)$ for some effective Cartier divisor $G$ on $Y$.
Then we have
\begin{equation}\label{STFormula}
\tau(\omega_X, \mathfrak{a}^t) =\sum_{e\geq 0}  \Tr^{e}\left( F^e_* \pi_*\OO_{Y}( \lceil K_{Y} - p^et G\rceil - D)   \right) 
\end{equation}
for any sufficiently large divisor $D$, and  $e>e_0$ depending on $D$.  
\end{Prop}

\begin{Rmk}
    Schwede and Tucker proved this elegant formula \eqref{STFormula} to compute test ideals of triples, we can get our ``test module"-version just by taking $\Delta = -K_X$ in their formula. 
\end{Rmk}

\section{Main Results}\label{section-applications of $F$-rationality}

\begin{Asm}\label{Char-p-Assumption}
Throughout this section $(R,\maxm)$ is an $F$-finite, reduced ring of dimension $d\geq 2$ and $\mathfrak{a}\subseteq R$ an ideal with positive height. We set $\rS=R[\mathfrak{a}t]$ to be the Rees algebra of $\mathfrak{a}$ with homogeneous maximal ideal $\maxm_S =  \maxm \oplus \mathfrak{a} \oplus \mathfrak{a}^2 \oplus \dots$ and $T=R[\mathfrak{a}t,t^{-1}]$ the extended Rees algebra of $\mathfrak{a}$ with homogeneous maximal ideal $\maxm_T$.
\end{Asm}

\begin{Th}\label{Th-test module decompositions}
  Under Assumption \ref{Char-p-Assumption}, for any $\lambda \geq 0$, the test submodule of ($S, (\fa.S)^\lambda)$ is given by
\[\tau(\omega_{S}, (\mathfrak{a}.S)^{\lambda}) = \bigoplus_{n\geq 1} \tau (\omega_R, \mathfrak{a}^{n+\lambda})\] and the test submodule of $(T, (t^{-1})^\lambda)$  is given by
\[ [\tau(\omega_{T}, (t^{-1})^\lambda)]_{\geq 0} = \bigoplus_{n \geq 0} \tau (\omega_R, \mathfrak{a}^{n+\lambda}).\]
 In particular, $ [\tau(\omega_{T}, (t^{-1})^\lambda)]_{\geq 1} = \tau(\omega_{S}, (\mathfrak{a}.S)^{\lambda}) $ and $[\tau(\omega_{T}, (t^{-1})^\lambda)]_0 = \tau(\omega_R , \mathfrak{a}^\lambda)$.    
\end{Th}

\begin{proof}

Since the test module does not change under normalization \cite[Chapter 5, Section 4, Exercises]{FBook}
we may replace $S$ and $T$ by their normalization. Set $X=\Spec R$, $Y=\Proj S$ and $\pi: Y \to X$ the blowup along $\mathfrak{a}$ so that $\mathfrak{a} \OO_Y= \OO_Y(-G)$. Note that under this notation we have that $\omega_S=\bigoplus_{n>0}H^0(Y,\omega_Y(-nG))$, by Remark \ref{RemarkSdS}.

First proving the statement about the test module of $S$:
We begin by choosing a test element $d\in S$. Then we have that
     
\begin{align*}
    \tau(\omega_S, (\mathfrak{a}S)^\lambda)&= \sum_{e\geq 0} \Tr_S^e(F^e_*(d \cdot (\mathfrak{a}S)^{\lceil\lambda p^e \rceil}\omega_S))\\
    &=\sum_{e\geq 0} \Tr_S^e(F^e_*d \cdot (\mathfrak{a}S)^{\lceil\lambda p^e \rceil}\bigoplus_{m>0}H^0(Y,\omega_Y(-mG)))\\
    &=\sum_{e\geq 0} \Tr_S^e(F^e_*\bigoplus_{m>0}H^0(Y,\omega_Y(-mG-\lceil\lambda p^e \rceil G-D)))\\ 
    &=\bigoplus_{n>0} \sum_{e\geq 0}  \Tr_R^e(F^e_*H^0(Y,\omega_Y(-n.p^e.G-\lceil\lambda p^e \rceil G-D))) \
\end{align*}
where $d\cdot \OO_Y= \OO_Y(-D)$ and the last equality follows from Proposition \ref{Prop-traceonrees}. Note that, apriori, we just have $\subseteq$ containment from 2nd to 3rd equality. However, as $(-G)$ is ample, we get equality for $e\gg 0$. Alternately, we can always enlarge the test element to get the reverse inclusion. Rephrasing in sheaf theoretic terms, the prior expression corresponds to the sheaf $\bigoplus_{n>0} \sum_{e\geq 0} \Tr_X
^e( \pi_* F^e_* \omega_Y(-np^eG-\lceil\lambda p^e \rceil G-D)))$. So by Proposition \ref{Prop-tau from blowup} (after again possibly enlarging $D$), we have that  
\begin{align*}
    \bigoplus_{n>0}\sum_{e\geq 0}  \Tr_X^e( \pi_* F^e_* \omega_Y(-np^eG-\lceil\lambda p^e \rceil G-D)))&=\bigoplus_{n>0} \sum_{e\geq 0} \Tr_X^e( \pi_* F^e_* \omega_Y(\lceil-(\lambda+n) p^e \rceil G-D)))\\
    &=\bigoplus_{n>0} \tau(\omega_X, \mathfrak{a}^{\lambda+n}),
\end{align*}
thus proving the claim.

Now addressing the claim on $T$: In positive degrees the claim is immediate from the fact that $(t^{-1})$ and $\mathfrak{a}T$ in non-negative degrees and that by Theorem \ref{Th-compareomegaSandT} we have that $\omega_S$ and $\omega_T$ agree in positive degrees.

For $k \geq 0$ and $e\gg0$, by Remark \ref{inclusionofgradedpieces}
\begin{align*}
    \left[ \tau(\omega_T, (t^{-1})^\lambda) \right]_{k} &= \sum_{e\geq 0} T^e_R \left( F^e_* \left( \left[ d\cdot  (t^{-1})^{\lceil\lambda p^e\rceil} \cdot \omega_T  \right]_{kp^e} \right) \right)\\
    &= \sum_{e\geq 0} T^e_R ( F^e_* ( d\cdot [\omega_T]_{\lceil\lambda p^e\rceil +kp^e}) ) \\  
\end{align*}

Now if either $\lambda>0$ or $k>0$,  we can choose $e$ large enough that \[[\omega_T]_{\lceil\lambda p^e\rceil +kp^e} =  [\omega_S]_{\lceil\lambda p^e\rceil +kp^e}\]

and hence with the same argument as above
\[\left[ \tau(\omega_T, (t^{-1})^\lambda) \right]_{k}= \tau(\omega_R,\mathfrak{a}^{\lambda+k}).\]
In particular, when $\lambda>0$, we have that \[\left[ \tau(\omega_T, (t^{-1})^\lambda) \right]_{0}= \tau(\omega_R,\mathfrak{a}^{\lambda}).\]


Finally the case of $\lambda=0$, $k=0$. Pick a test element $ct^{-r}$ for $r\gg0$. Then, for $D = div(c)$, by using Remark \ref{inclusionofgradedpieces},

\begin{align*}
\left[ \tau(\omega_T) \right]_0 
&= \left[ \sum_{e \geq 0} \Tr^e_T \left( F^e_*{(\omega_T \cdot c t^{-r})} \right) \right]_0 \\
&= \sum_{e \geq 0} \Tr^e_R \left( F^e_*[\omega_T]_{r} \cdot c \right) \\
&= \sum_{e \geq 0} \Tr^e_R \left( F^e_*{\pi_* \cO_Y(K_Y - rG - D)} \right) \\
&\subseteq \sum_{e \geq 0} \Tr^e_R \left( F^e_*{\pi_* \cO_Y(K_Y - D)} \right) \\
&= \tau(\omega_R).
\end{align*}

giving containment on one side. For the opposite direction, note that there exists a Cartier divisor $\Gamma>0$ such that:

\[ \pi_* \cO_Y(K_Y - rG - D) \supseteq \pi_* \cO_Y(K_Y - D - \Gamma)\]
 giving, 
\[ \sum_{e \geq 0} \Tr^e_R \left( F^e_*{\pi^e_* \cO_Y(K_Y - rG - D)} \right) \supseteq \sum_{e \geq 0} \Tr^e_R \left( F^e_*{\pi^e_* \cO_Y(K_Y - D - \Gamma)} \right) =  \tau(\omega_R).\]

(see \cite[Proof of Lemma 3.5]{Schwede-Tucker2014Testidealsofnon-principalideals:computationsjumpingnumbersalterationsanddivisiontheorems} for example)

Here, "$D+\Gamma$" can be interpreted as enlarging the test element.

This establishes the equality $\left[\tau(\omega_T)\right]_0 = \tau(\omega_R)$.
\end{proof}

We can also compute the negative degrees of the test module of $T$ as follows. First we need a lemma.

\begin{Rmk}\label{LambdaLimit}
    Let \(\lambda \geq 0\) be a real number and \(p\) be a prime. Then,
\[
\lim_{e \to \infty} \frac{\lceil p^e \lambda \rceil}{p^e} = \lambda
\]
\end{Rmk}

\begin{Th}\label{Th-testmoduleExtReesNeg} For $n \leq -1$
    \[[\tau(\omega_T, (t^{-1})^\lambda)]_n = 
    \begin{cases} 
\tau(\omega_R) , & \forall n \leq -\lambda \\
\tau(\omega_X, \mathfrak{a}^{\lambda+n}), & -\lambda < n \leq -1
\end{cases}
    \]
\end{Th}

\begin{proof}
    Let $d \in $R$ = [T]_0$ be a test element for both $R$ and $T$. Then rearranging sums and applying Proposition \ref{Prop-traceonrees} and Remark \ref{inclusionofgradedpieces},
    \begin{align*}
        [\tau(\omega_T, (t^{-1})^\lambda)]_n &= \left[ \sum_{e \geq 0} \Tr_T F_*^e \omega_T (t^{-1})^{\lceil p^e \lambda \rceil} d  \right]_n\\
        &=  \sum_{e \geq 0} \Tr_R [F_*^e  \omega_T (t^{-1})^{\lceil p^e \lambda \rceil} d  ]_{np^e}\\
        &=  \sum_{e \geq 0} \Tr_R [F_*^e  \omega_T ]_{np^e + \lceil p^e \lambda \rceil}d.
    \end{align*} 

  By the Lemma \ref{LambdaLimit}, we will consider the following 3 cases.  
    
   \textbf{Case-1:} $\mathbf{\lambda +n < 0.}$\\
   By Proposition \ref{Prop-negativeeonghdegrees} we have that $[\omega_T]_k$ agrees with $\omega_R$ for all $k\leq k_0$. We may then replace $\sum_{e \geq 0}$ in the above sum with $\sum_{e \geq e_0}$ so as to make $np^e + \lceil p^e \lambda \rceil \leq k_0$. Thus,
    \begin{align*}
        \phantom{\tau(\omega_T, (t^{-1})^\lambda)]_n} &=  \sum_{e \geq e_0} \Tr_R F_*^e \omega_R d\\
        &= \tau(\omega_R)
    \end{align*}

     \textbf{Case-2:} $\mathbf{\lambda +n > 0.}$\\
     Then by Theorem \ref{Th-compareomegaSandT}, $[\omega_T ]_{np^e + \lceil p^e \lambda \rceil} = [\omega_S ]_{np^e + \lceil p^e \lambda \rceil} = H^0(Y,\omega_Y(-n.p^e.G-\lceil\lambda p^e \rceil G)) $ by Remark \ref{RemarkSdS}. So, we get 
     \begin{align*} 
     =\sum_{e\geq 0}  \Tr_R^e(F^e_*H^0(Y,\omega_Y(-np^eG-\lceil\lambda p^e \rceil G-D))) = \tau(\omega_X, \mathfrak{a}^{\lambda+n}). \\
     \end{align*}

     \textbf{Case-3:} $\mathbf{\lambda +n = 0.}$\\
     Then we basically repeat what we did for the case $\lambda = k = 0$ in the previous Theorem \ref{Th-test module decompositions}.
Pick a test element $ct^{-r}$ for $r\gg0$. Then,

\begin{align*}
&= \sum_{e \geq 0} \Tr^e_R \left( F^e_*[\omega_T]_{r} \cdot c \right) \\
 &=\sum_{e\geq 0}  \Tr_R^e(F^e_*H^0(Y,\omega_Y(-rG-D))) \\
&= \sum_{e \geq 0} \Tr^e_R \left( F^e_*{\pi_* \cO_Y(K_Y - rG - D)} \right) \\
&\subseteq \sum_{e \geq 0} \Tr^e_R \left( F^e_*{\pi_* \cO_Y(K_Y - D)} \right) \\
&= \tau(\omega_R).
\end{align*}

For the opposite inclusion, note that there exists a Cartier divisor $\Gamma>0$ such that:

\[ \sum_{e \geq 0} \Tr^e_R \left( F^e_*{\pi^e_* \cO_Y(K_Y - rG - D)} \right) \supseteq \sum_{e \geq 0} \Tr^e_R \left( F^e_*{\pi^e_* \cO_Y(K_Y - D - \Gamma)} \right) = \tau(\omega_R).\]

So, the proof is complete.
     
\end{proof}

\begin{Rmk}
    
    Similar decomposition for the Rees Algebra $S$ is obtained by Hara and Yoshida in \cite[Theorem 5.1]{HaraYoshida} under the assumption that $\fa$ is $\maxm$-primary and Proj $S$ is $F$-rational. Kummini and Kotal generalized Hara-Yoshida's  decomposition in \cite{KotalKummini} under the assumption that $S$ is Cohen-Macaulay and the $\maxm$-primary ideal $\fa$ has a reduction generated by a system of parameters.
\end{Rmk}

\begin{Rmk}\label{testmodulefor triples}
    Essentially following the same argument and the fact that for a principal divisor $D = (f)$, $\tau(R, \Delta+D, \mathfrak{a}^t)=f.\tau(R, \Delta, \mathfrak{a}^t)$, one can get similar decompositions for the test module of triples, just by keeping track of grading. But we don't include it here for the sake of brevity.
\end{Rmk}

\begin{Rmk}\label{testidealdecomposition}
    By using Lemma \ref{Relatingtestmoduleandideal} and Remark \ref{testmodulefor triples} one could also get similar decomposition for test ideals (for triples) as well. In particular, when $R$ is $\bbQ$-Gorenstein, we will have $[\tau(T,(t^{-1})^\lambda)]_0 = \tau(R , \fa^\lambda) $.
\end{Rmk}

\section{Equivalence of $F$-rationality of Rees and extended Rees Algebras}
Throughout this section we maintain the notation of the previous section, that is: We let $(R,\maxm)$ be an $F$-finite, reduced ring of dimension $d\geq 2$ and $\mathfrak{a}\subseteq R$ an ideal with positive height. We set $\rS=R[\mathfrak{a}t]$ to be the Rees algebra of $\mathfrak{a}$ with homogeneous maximal ideal $\maxm_S =  \maxm \oplus \mathfrak{a} \oplus \mathfrak{a}^2 \oplus \dots$ and $T=R[\mathfrak{a}t,t^{-1}]$ the extended Rees algebra of $\mathfrak{a}$ with homogeneous maximal ideal $\maxm_T$.

For an $\mathfrak{m}$-primary ideal $\mathfrak{a}$ it was proven by Hara, Watanabe, and Yoshida that when $R$ and $S$ are $F$-rational then so is $T$ \cite[Theorem 4.2]{HWY-02F-Rational}, and conjectured the converse \cite[Conjecture 4.1]{HWY-02F-Rational}  which was proven by Koley and Kummini who also provided a new proof of the forward direction \cite{Koley-Kummini21}. Using Theorem \ref{Th-test module decompositions} and Theorem \ref{Th-testmoduleExtReesNeg}, we answer Hara, Watanabe and, Yoshida's Conjecture as well as give a different proof of the forward direction to any ideal (not necessarily $\mathfrak{m}$-primary).

\begin{Th}\label{Th-equivofFrationality}
    Under Assumption \ref{Char-p-Assumption}, we have the following implications:
    \begin{enumerate}
        \item $R$ and $S$ are both $F$-rational implies $T$ is $F$-rational.
        \item $T$ is $F$-rational implies $R$ and $S$ are $F$-rational.
    \end{enumerate}
\end{Th}

\begin{proof}
    First we prove the equivalence about Cohen-Macaulayness. Note from \cite[Proposition 1.1]{Huneke-82}  it follows that $R$ and $S$ are Cohen-Macaulay implies $G= T/(t^{-1})$ is Cohen-Macaulay and hence $T$ is Cohen-Macaulay as well. To prove the other direction, note the following:
    \begin{itemize}
        \item $T$ is Cohen-Macaulay \iff $G$ is Cohen-Macaulay.
        \item $T$ is $F$-rational implies $ T_{t^{-1}} = R[t, t^{-1}]$ is $F$-rational implies $R$ is $F$-rational and hence, $R$ is pseudo-rational by \cite{SmithF-rationalPseudorational}. 
        \item $T$ being Cohen-Macaulay and $R$ being pseudo-rational implies, from Theorem (5) of \cite{Lipman1994CM} that $S$ has to be Cohen-Macaulay as well.
    \end{itemize}

    Hence, we need only compare the test modules of $R$, $S$, and $T$. With this in mind (2) is immediate from Theorem \ref{Th-test module decompositions} by putting $\lambda = 0$.
    
    On the other hand (1) is clear from Theorem \ref{Th-test module decompositions}, Theorem \ref{Th-testmoduleExtReesNeg} and Proposition \ref{Prop-canonicalnegativedegrees}.
\end{proof}

\begin{Rmk}\label{deformation}
    If $G = T/(t^{-1})$ is $F$-rational, then we get $T$ is $F$-rational and hence $S$ is $F$-rational  \cite[see page 182]{HWY-02F-Rational}. Similar results were known for Cohen-Macaulay and Gorenstein property by extensive works of Goto, Huneke, Ikeda, Lipman, Shimoda, Vi\'et etc, see \cite{GotoShimodaCM}, \cite{iai2024characterizationsgorensteinreesalgebras}, \cite{Goto-Nishida-Book}, \cite{Ikeda-Gor}, \cite{VietCM}, \cite{WhenCM}, \cite{WhenGor} \cite{Huneke-82} and, \cite{Lipman1994CM} for example.
\end{Rmk}

\section{Other Applications}\label{sec-other applications}

\begin{App}\label{Macaulay2}
    In a separate paper, Karl Schwede and the authors will apply Theorem \ref{Th-test module decompositions} (specifically the equality $[\tau(\omega_{T}, (t^{-1})^\lambda)]_0 = \tau(\omega_R , \mathfrak{a}^\lambda)$) to compute test ideals of non-principal ideals of a $\bbQ$-Gorenstein ring in Macaulay2. This would confirm Schwede and Tucker's prediction in \cite[Remark 5.2]{Schwede-Tucker2014Testidealsofnon-principalideals:computationsjumpingnumbersalterationsanddivisiontheorems}. Note that for principal ideals, there is already a testideal package available. 
\end{App}

\begin{App}\label{F-Jump[ing]}
    If we assume the principal ideal case \cite[Chapter-4, Section 6.3]{FBook}, then we can use Remark \ref{testidealdecomposition} to give a quick proof of discreteness and rationality of $F$-jumping numbers for non-principal ideals when $R$ is $\bbQ$-Gorenstein, recovering \cite[Theorem B]{Schwede-Tucker2014Testidealsofnon-principalideals:computationsjumpingnumbersalterationsanddivisiontheorems}. Also note that, it directly follows from Remark \ref{testidealdecomposition} that if $\lambda$ is an $F$-jumping number for $(R, \fa)$ then so is $\lambda + n$ for infinitely many $n \geq 0$, recovering \cite[Proposition 1.12]{ELSV}.
\end{App}

\begin{App}
If $R$ is essentially finite type over a field of characteristic 0, by the works of Hara \cite{Hara}, Mehta-Srinivas \cite{MehtaSrinivas} and Smith \cite{SmithF-rationalPseudorational}, our Theorem \ref{Th-equivofFrationality} immediately recovers the corresponding statement for rational singularity, giving an alternate proof of \cite[Theorem 3.2]{HWY-F-Regular}.
\end{App}

There are many other potentially interesting future directions to apply our methods to simplify and generalize various positive charactaristic results like \cite{Mustata-Zhang}, \cite{BS-Mustata}, \cite{Sato} and, \cite{smirnov2023theoryfrationalsignature}. However, we will leave these as future projects.

 \bibliographystyle{amsalpha}
  \bibliography{refs}

\end{document}